\documentclass{amse-new}

\graphicspath{{figures/}}

\numberwithin{equation}{section} 

\begin{document}

 \PageNum{1}
 \Volume{Peking University}{}{}{}
 \OnlineTime{School of Mathematical Sciences}
 \DOI{Minzheng Li}
 \EditorNote{Received x x, 201x, accepted x x, 201x}

\abovedisplayskip 6pt plus 2pt minus 2pt \belowdisplayskip 6pt
plus 2pt minus 2pt
\def\vsp{\vspace{1mm}}
\def\th#1{\vspace{1mm}\noindent{\bf #1}\quad}
\def\proof{\vspace{1mm}\noindent{\it Proof}\quad}
\def\no{\nonumber}
\newenvironment{prof}[1][Proof]{\noindent\textit{#1}\quad }
{\hfill $\Box$\vspace{0.7mm}}
\def\q{\quad} \def\qq{\qquad}
\allowdisplaybreaks[4]


\AuthorMark{Minzheng Li}                             

\TitleMark{The Continuity Equation Contains Non-Stochastic Motion}  

\title{The Continuity Equation Contains Non-Stochastic Motion        
\footnote{Supported by \ldots (Grant No. \ldots)}}                  

\author{Minzheng Li}             
    {Address: School of Mathematical Sciences, Peking University \\
    E-mail 1801110057@pku.edu.cn \,$ $ }

\maketitle%

\Abstract{The scientific question resolved by this paper is that the continuity equation appears as an equivalent language of the system of first-order linear ODE. The main result characterizes the fact that the continuity equation contains non-stochastic motion; a stochastic motion addressed by the continuity equation surely drops its stochastic part in a probabilistic indistinguishable manner.}      


\MRSubClass{34Axx,62A01}      

\def\d{\mathrm{d}}
\def\v{\mathbf{v}}
\def\f{\mathbf{f}}
\def\P{\mathbb{P}}
\def\E{\mathbb{E}}
\def\I{\mathrm{I}}
\def\scrX{\mathscr{X}}
\def\scrY{\mathscr{Y}}
\def\scrZ{\mathscr{Z}}
\def\R{\mathbb{R}}
\def\var{\mathrm{var}}

\section{Introduction}
The Reynolds transport theory is presented in many well-written fluid mechanics textbooks. It is a phenomenon that particle motion in first-order ODE includes a statistical behavior, named as the continuity equation. However, whether the continuity equation contains or partially contains information on particle motion is to be questioned. After proving that the continuity equation is unique to the velocity field it contains, the question is specified as whether the continuity equation guarantees non-stochastic motion. In this paper, we pursue reversed proposition of the Reynolds transport theory, which together with Reynolds transport theory constitutes the fact that the continuity equation appears as the third equivalent mathematical language for the linear first-order ODE. The essence of the proof of our result is that the continuity equation has characterized arbitrarily small variance on its own, ruling out the possibility of existence of stochastic motion. The findings provide a new tool in the study of mathematical ordinary differential equations and the study of descriptive statistics and probabilistic analysis.

The system of first-order linear ODE\begin{equation}
\label{Abstract_ODE}\frac{\d}{\d t}x(t) = \v(x(t),t)
\end{equation}
is the fundamental form for ODE theories. Here $ x(t)\in \R^p $ is motion in $p$-dimensional Euclidean space. $ \v\in \R^p\times\R\to \R^p $ is time-dependent velocity field. In the cases that $ \v $ explicitly contains time $ t $, we address the system as non-autonomous. In the cases that $ \v $ does not contain $ t $, the system is autonomous. Continuous-time Newton's method is an example of the autonomous system:
$$
\dot{x} = - J^{-1}(x)F(x)
$$
where $ F $ is the target root-searching function or the gradient function of the target optimization function. $$
J = \mathrm{grad}(F) = \triangledown F
$$is the gradient function of the target root-searching function $ F $, or the Hesse matrix of the target optimization function. Note that when $ p = 1$ $ F $ and $ J $ are both numerical functions with one variable, and when $ p \geqslant 2$, $ F $ is a multi-variable vector-valued function and $ J $ is $ p\times p $ matrix. In the common basic theory of ODE, linear ODE of order higher than one can revert to the system of first-order linear ODE by techniques of change of variables. This technique is typically applied in physics when we write the second-order derivatives of position into two independent variables position and momentum.

Another analytic form of \eqref{Abstract_ODE} is \begin{equation}
\label{Abstract_ODE_2}
x(t) = x(0) + \int_0^t \v(x(s),s)\d s
\end{equation}
It is clear to see the equivalence of \eqref{Abstract_ODE} and \eqref{Abstract_ODE_2}: \eqref{Abstract_ODE} is \eqref{Abstract_ODE_2} to take derivative of time, and \eqref{Abstract_ODE_2} is \eqref{Abstract_ODE}to take the integral of time. Note \eqref{Abstract_ODE_2} is to take the integral of a vector-valued function. The importance of the analytic form \eqref{Abstract_ODE_2} is of the theorem of existence and uniqueness of solutions of ODE, also known as Picard's iteration theorem. The proof, as its name suggests, is to exploit \eqref{Abstract_ODE_2} into a structure like contraction mapping, with the condition of Lipschitz continuity of the velocity field. In our development of theory, the Lipschitz continuity condition of the velocity $ \v $ is presumptively provided, thus frees us of the possibility of existence of multi-solutions of the particle movement.

Is there the third additional equivalent form of ODE \eqref{Abstract_ODE} and \eqref{Abstract_ODE_2} to describe particle motion in the velocity field? This is exactly the question to be answered in this paper. We declare yes and the continuity equation comprehending all its possible initial values is the third equivalent language for the general system of linear first-order ODE. The main result is that the continuity equation contains non-stochastic motion.

We start our statement with the conception of fluid in physics. In physics, the macro fluid is depicted by the density function. Physicists believe that fluid is continuous media, though it is composed of a large number of discrete particles. How large is the number, consider 1 mol in chemistry, and it is around $6.02\times 10^{23}$ number of particles. Such a quantity of particles amounts to one spoon in laboratories. Fluids are compressible fluids and incompressible fluids. In our everyday life air is compressible and water is incompressible. Consider an ensemble of $ N \gg 1 $ $ p $-dimensional particles $ x^1(t), \cdots, x^N(t) $. $ N $ is extremely large such as $ N = $ 1000 mol or 1000000 mol. The ensemble of $ x^1(t), \cdots, x^N(t) $ is considered as a compressible macro continuous fluid. Define $$
\int_B \rho(x;t)\d x = \frac{1}{N}\sum\limits_{j=1}^N \I (x^j(t)\in B)
$$
$ B \subset \R^p $ is subset in $p $-dimensional Euclidean space and $ \I $ is the indicator function. $ \rho $ is the statistical density function of the ensemble $ x_1(t), \cdots, x_N(t) $, satisfying $1$ as its total integration. $ \rho $ is viewed as continuous because $ N $ is extremely large. Here the concept of $ p $-dimensional particles is purely mathematical. In the real physical world, all the experiments are conducted in three-dimensional space with half-dimensional time. There is no proof of the existence of dimension higher than three in the real world. We demand $ x^1(t), \cdots, x^N(t) $ have unequal initial values $ x^1(0), \cdots, x^N(0) $, thus according to the theorem of existence and uniqueness of solutions of ODE or the Picard's iteration theorem, $ x^1(t), \cdots, x^N(t) $ do not intersect at any time $ t $.

Reynolds transport theory states that when $ N \gg 1 $ the statistical density of $ x^1(t), \cdots, x^N(t) $, $ \rho(x;t) $ satisfies the continuity equation
\begin{equation}
\label{Abstract_Continuity_Eq2}
\frac{\partial\rho(x,t)}{\partial t} + \mathrm{div}_x(\mathbf{v}(x,t)\rho(x,t))=0
\end{equation}
or written as
\begin{equation}
\label{Abstract_Continuity_Eq}
\frac{\partial\rho(x,t)}{\partial t} + \triangledown_x\cdot(\v(x,t)\rho(x,t))=0
\end{equation}
where $ \mathrm{div}_x $ and $ \triangledown_x\cdot $ denotes taking divergence of the spatial variable $ x $. $ \v $ is the time-dependent velocity field of the ODE. The theorem as special cases holds for time-independent velocity fields. $ x(t) $ denotes the particle, $ x $ denotes location; they are irrelevant symbols. Note that $ \rho $ is not differentiable by its definition. Physicists tackle this issue by ignoring it; because the requirement of the theory is $ N $ being sufficiently large, such as 1000 mol or $ 10^6 $ mol, where $ x^1(t), \cdots, x^N(t) $ as a whole is seen as continuous media. The smooth function $ \rho $ describes the physical world adequately, and the differentiability of $ \rho $ is perfect for mathematical analysis. The mathematicians tackle this issue by assuming a smooth version $ \tilde{\rho} $ similar to $ \rho $, such as the convolution of $ \rho $ with a smooth mollifier. We follow the ideology of physicists and ignore this issue.

What we prove is the contrary proposition of Reynolds transport theory: For any initial value $ \xi(0) $, particles satisfy the stochastic motion of 
\begin{equation}
\label{Abstract_Stochastic_Motion}
\d \xi(t) = \v^*(\xi(t),t)\d t + \sigma^*(\xi(t),t)\d W(t)
\end{equation}
where $ W(t) $ is arbitrary stochastic process. Writing in full probabilistic language this is 
 $$
\d \xi(t;\omega) = \v^*(\xi(t;\omega),t)\d t + \sigma^*(\xi(t;\omega),t)\d W(t;\omega)
$$
The statistical density $ \rho(x;t) $ is defined as $$
\int_B \rho(x;t)\d x = \frac{1}{N}\sum\limits_{j=1}^N \I (\xi^j(t)\in B)
$$
$ N $ is extremely large such as $ 10^6 $ mol so that $ \rho $ is conceived to be smooth. Give $ \xi(0) $ distribution $ \rho(\cdot;0) $. If for any of initial distributions of $ \xi^1(0),\cdots,\xi^N (0)$ characterised by the statistical density $ \rho(x;0) $, $\rho(x,t)$ satisfies the continuity equation \eqref{Abstract_Continuity_Eq2} or written as \eqref{Abstract_Continuity_Eq}, then we have$$
\v^*(x,t)\equiv \v(x,t)
$$
and $$
\sigma^*(x,t) = 0
$$
In probabilistic language, it is guaranteed that particles do non-stochastic motion in the velocity field $ \v $, which is uniquely characterized by the continuity equation \eqref{Abstract_Continuity_Eq}; and \eqref{Abstract_Stochastic_Motion} and \eqref{Abstract_ODE} are probabilistically indistinguishable processes. We have introduced that \eqref{Abstract_ODE_2} and \eqref{Abstract_ODE} are equivalent languages, and now the continuity equation \eqref{Abstract_Continuity_Eq2} or written as \eqref{Abstract_Continuity_Eq} is the third equivalent language of the ODE \eqref{Abstract_ODE}. In the case that $ W(t) $ is Brownian motion $ B(t) $, the result is characterized by the classical Fokker-Planck equation in probability. But the result we provide is the case $ W(t) $ being any stochastic process such as $$
W(t) = \|B(t)\|^2 B(t)
$$ or $$
W(t) = \|B(t)\|^4 B(t)
$$
$ W(t) $ may not have a variance of growth proportional to time elapse. Our result is the contrary proposition of the Reynolds transport theory: the particle motion is guaranteed by its statistical nature.

The exemplification of the importance of equivalent languages should give to the discovery of the equivalence of the Schr\"{o}dinger equation and matrix mechanics. This unified the initial theoretic physics of the micro world and later gave rise to a method of quantization.

Equivalence is always the pursuits of mathematicians. Mathematical researchers develop good theories by stating several definitions or results are equivalent, and these would usually be the most classic results to write in textbooks. The system of first-order linear ordinary differential equations is equivalent to its statistical description the continuity equation, in the condition that the continuity equation comprehends all its possible initial distributions to start with. We hope that this equivalence provides a new possible approach to the study and the proofs of ODE theories.

In the confrontation with machine learning and artificial intelligence, the effectiveness of computerized models is expected to be studied by statisticians. On the other hand in the confrontation with real-world data, applying ideas for experiment design and causal inference is emphasized. And the writer believes that equivalent languages of ODE and the continuity equation would be important tools in the development of mathematical statistics.

The author expects that the statistical value of our results is higher than the mathematical. In most scientific subjects, description is the first step of the study. In psychology, psychologists need to describe behaviors and then predict behaviors. In physics, physicists need to describe the world and then predict the outcomes of experiments. It is the same for statistics. To study the phenomena of how operators behave in machine learning and artificial intelligence, the development of equivalent languages to address phenomena is key.

Macro physics and particle physics are two systems in the study of physics. Theories such as electromagnetics and fluid dynamics study macro statistical phenomena, and theories such as classical dynamics and quantum mechanics study the existence and movement of particles. Fluid dynamics do not consider non-physic scenarios like "a small iron ball flowing in the air stream". The so-called continuity equation in fluid dynamics is based on the continuity assumption, that the media is continuous and almost uniformly shaped. Our theory develops on a velocity field; however, in classic Navier-Stokes equations, the velocity functions are not velocity fields for particles but statistical descriptions of the macro movement. In this text, "a compressed small 'iron' ball flowing in air stream by the velocity field" is what we do consider, a violation of practical physics but key in mathematics.

There are of course many results in physics trying to link the particle with the macro. Liouville's theorem in statistical physics targets Hamilton-conserving particles. The Boltzmann equations, usually written as $$
\frac{\partial}{\partial t}\rho(x,v,t) + v\cdot \frac{\partial}{\partial x}\rho(x,v,t) + F\cdot \frac{\partial}{\partial v}\rho(x,v,t) = -\tau (\rho -\rho^{(eq)})
$$ $ \rho^{(eq)} $ being the local equilibrium state, appears in many modeling problems. DFT(Density Functional Theory)-based quantum chemistry enables scientists to simulate the structure of crystals or the molecular dynamics of liquids from the quantum level. \cite{PhysRev.140.A1133} proposed the well-known Kohn-Sham equations, regarded as a cornerstone of DFT, but to take too high computational costs. \cite{osti_1542041} takes a different approach by learning potential functions into neural networks. By the time the text is written no connection with these previous results has been found.

\section{Notation convention}
In the whole text, $ \mathrm{div} = \triangledown\cdot $ and $ \triangledown $ denote taking divergence and taking gradient of a vector function or a function of its location part $ x $ but not including its time part $ t $. $ x $ is $ p $-dimensional and time $ t $ is one-dimensional. $ \rho(x;t) $ is exactly the same with $ \rho(x,t) $ but written simply to emphasise time $ t $ is time parameter. $ dt $ is essentially exactly the same as $ \d t $. We write $ dt $ simply to emphasize it is a very instant time forward considered, and $ \d t $ emphasizes the mathematical calculation of differentiation and integration. In most ODE language $ \v $ is written as $ \v(t,x) $, but to comply with the writing of $ \rho(x,t) $, the language seems to be better to write as $ \v(x,t) $. In the whole text $ \v(t,x) $ and $ \v(x,t) $ are exactly the same but different language styles to write. Note that $$
\triangledown\cdot\f = \frac{\partial\f_1}{\partial x_1} + \cdots +\frac{\partial \f_p}{\partial x_p}
$$
is not equal to $$
\f\cdot\triangledown = \f_1\frac{\partial}{\partial x_1} + \cdots + \f_p\frac{\partial }{\partial x_p}
$$
That is to say in $ \triangledown_x\cdot $, $ \cdot $ does not satisfy the definition of inner product. However, since the writing of $ \triangledown_x\cdot $ is so popular in physics that we follow this writing and deprive ourselves sometimes of writing $ \mathrm{div}_x $, the symbol mostly in mathematics. We abuse the symbol slightly: $ x(t) $ with superscript $ 1,\cdots, N $ denotes $ p$-dimensional particles, whereas $ x $ denotes $ p $-dimensional position, its components $ x = (x_1,\cdots,x_p) $ are presented by subscripts. $ x(t) $ and $ x $ are irrelevant. $ x $ with additional $ (t) $ $ x(t) $ denotes particles and its motion with time $ t $, whereas $ x $ is simply location. $ \v_j $ is the $ j $-th component of $ \v(x) = \begin{bmatrix}
\v_1(x_1, \cdots, x_p)\\
\vdots\\
\v_p(x_1, \cdots, x_p)\\
\end{bmatrix} $.
$ \v_j(x_1, \cdots, x_p) $ is a scalar-valued function. Without specific explanation, the mathematical calculations in this text are to both scalar-valued functions and vector-valued functions, and even to matrix-valued functions. For example, taking gradient to a vector $ u $ $$
\triangledown_u u
$$
yields an identity matrix of the same dimension of $ u $. An operator $ (\partial_t+\v(t,u)\partial_u) $ is introduced, which abbreviates $(\partial_t+\v(t,u)\cdot\partial_u)$ to reduce space. Because the rationale of the proof is to set $ \Sigma $ the variance of $ \rho(\cdot;t) $ arbitrarily small and observe how the concentration of mass behaves, we may omit writing $ o(\Sigma) $ to reduce space, in the hope that readers pardon us for possible reduction of writing.

\section{Main Results}
We shall give our main results in this section and its proof in the next section, and later provide a discussion of the details in the proof.

\begin{theorem}\label{MainTheorem}
\textbf{First}, particles moving in the time-dependent velocity field $ \mathbf{v}(t,\cdot) \in C(\mathbb{R}\times\mathbb{R}^p,\mathbb{R}^p ) $ is written in ODE \begin{equation}
\label{ODE}
\frac{\d x(t)}{\d t} = \mathbf{v}(t,x(t))
\end{equation} Consider an extremely large amount of particles moving independently, charaterized by their statistical density $ \rho(x,t)\in C^1(\mathbb{R}^p\times\mathbb{R},\mathbb{R} ) $. We have that $ \rho(\cdot, \cdot) $ satisfies the continuity equation\begin{equation}
\label{TheContinuityEquation}
\frac{\partial}{\partial t}\rho(x;t) + \triangledown\cdot\left(\rho(x;t) \mathbf{v}(t,x)\right)=0
\end{equation} where $ \triangledown\cdot = \mathrm{div} $ is taken on location part $ \mathbb{R}^p $. This is the classic Reynolds transport theory.

\textbf{Second}, the velocity field $\v $ that the continuity equation \eqref{TheContinuityEquation} contains is unique. In standard mathematical language if $ \rho $ satisfies both $$
\frac{\partial}{\partial t}\rho(x;t) + \triangledown\cdot\left(\rho(x;t) \mathbf{v}(t,x)\right)=0
$$
and $$
\frac{\partial}{\partial t}\rho(x;t) + \triangledown\cdot\left(\rho(x;t) \mathbf{v}^*(t,x)\right)=0
$$
for its any initial value $ \rho(\cdot;0) $, then we have
$$
\mathbf{v}(\cdot,\cdot) \equiv \mathbf{v}^*(\cdot,\cdot) 
$$ 

\textbf{Third}, we are jumping our discussion from physics into pure mathematics. We are addressing what the continuity equation \eqref{TheContinuityEquation} acts like by setting $ \rho $ almost identical to Dirac delta function, note that $ \rho $ is normalized to one as its total integration by its definition in Introduction. This is the physical scenario to compress the whole air in the room to a small "iron" ball and observe where the small iron ball would move according to motion of \eqref{TheContinuityEquation}. Of course in physics neither we care nor we can do an experiment like this. The idea of compressing the whole air in the room to a small ball, which still obeys the continuity equation, is pure mathematical conjecture.

We state that with the statistical density $ \rho $ satisfying the PDE \begin{equation}\label{MainResult_eq2}
\frac{\partial}{\partial t}\rho(x;t) + \mathrm{div}\left(\rho(x;t) \mathbf{v}(x,t)\right)=0
\end{equation}
a particle motion along the direction of $ \mathbf{v}(t,\cdot) $ is determined. The small iron ball moves exactly in the velocity field $ \v $ for an instant time ahead. Of course $ \v $ is uniquely designated by the continuity equation, as we just declared in the second part of this theorem.

In the fourth and fifth parts of the theorem, we are stating, in complete probabilistic language, a contrary proposition of Reynolds transport theorem. \textbf{Fourth}, consider independent and identically distributed diffusion processes $ \xi(t) $ \begin{equation}
\label{WienerProcess}
\d\xi(t) = \mathbf{v}^*(t,\xi(t))\d t + \sigma^*(t,\xi(t))\d B(t)
\end{equation}
$ B(t) $ being Brownian motion. If their statistical density $ \rho $ satisfies \eqref{TheContinuityEquation} with its arbitrary initial value $ \rho(\cdot;0) $, then we have $$
\mathbf{v}^*(t,x) \equiv \mathbf{v}(t,x) 
$$and $$
\sigma^*(t,x) \equiv 0 
$$ 
$ \xi $ is probabilistically indistinguishable process with the trajectory of ODE of \eqref{ODE}. This is the result of the classic Fokker-Planck equation. Here $ \rho $ can be understood in two ways: $ \rho $ is the statistical density of i.i.d. samples of $ \xi $: $ \xi^{(1)},\cdots,\xi^{(N)}$, $ N =1000 $ mol; or the probability density of $ \xi $. Both understandings are OK, and the result is the same.

The fourth part of the theorem has constrained the stochastic part of motion $ \xi $ to be Brownian motion, whereas our result is this to be any possible stochastic motion. \textbf{Fifth}, to generalize our result to all stochastic processes, we need a requirement that  $ \v $ in the continuity equation is smooth $$
\v\in C^\infty(\R^p\times\R;\R^p)
$$
and the statistical density $ \rho $ and the solution of ODE \eqref{ODE} are analytic. Given that the statistical density $ \rho $ of i.i.d. stochastic processes $ \xi(t) $
\begin{equation}
\label{StochasticProcess_Mainresults}
\d\xi(t) = \mathbf{v}^*(t,\xi(t))\d t + \sigma^*(t,\xi(t))\d W(t)
\end{equation}
where $ W(t) $ is any mean-zero stochastic process such as $$
W(t) = \|B(t)\|^4 B(t)
$$ or $$
W(t) = \|B(t)\|^6 B(t)
$$
satisfies the continuity equation \eqref{TheContinuityEquation}, for any its initial value $ \rho(\cdot;0) $ the statistical density of $ \xi(0) $, then we have $$
\mathbf{v}^*(t,x) \equiv \mathbf{v}(t,x) 
$$and $$
\sigma^*(t,x) \equiv 0 
$$  
$ \xi $ is probabilistically indistinguishable process with exactly the same trajectory of ODE \eqref{ODE}.
\end{theorem}

\begin{remark}
When we state that "particles move in the velocity field $ \mathbf{v}(t,\cdot) \in C(\mathbb{R}\times\mathbb{R}^p,\mathbb{R}^p ) $" we presume the existence and uniqueness of the solution of \eqref{ODE} guaranteed. Please be noted that this may not happen for every ODE. Consider this initial value problem: 
$$
\left\{ \begin{aligned}
\frac{\d x}{\d t} = x^\frac{1}{2}\\
x(0) = 0
\end{aligned}
\right.
$$
Both $ x\equiv 0  $ and $ x = \frac{1}{4}t^2$ are solutions. Observing the two curves, Yes they share the same velocity at time 0, but they choose different velocities at the very instant moment they leave 0. To avoid troubles like this, we set out the uniqueness of the solution of ODE in our development of theory. 
\end{remark}

\begin{remark}
The first part of Theorem \ref{MainTheorem}, particle nature containing statistical nature, coincides with Reynolds transport theorem. The fourth part of theorem is the result of the Fokker-Planck equation. The rest parts of the proof are our work.
\end{remark}

\begin{remark}
We presume the uniform boundedness of the velocity field $ \v $: $$
\sup_{x\in\R^p}|\v(x)| < \infty
$$
This presumption is reasonable. $ \rho(x,t) $ is defined on $ \R^p\times\R $, $ \R $ is the time dimension, and $ \R^p $ is the location space; however in most real-world settings such as numerical solution of Newton's method, we are considering a compact set $ E\subset\R^p $. We have$$
\rho(x,t) = 0,x\notin E
$$
or
$$ \v(x,t) = 0,x\notin E
$$
satisfied. Typically in our proof, we believe our theory happens in the compact set $ E\subset \R^p $ and an instant time interval $ [t,t+dt]\subset[0,T] $. Compact sets such as $ E \times [0,T] $ are the time-location of interest.
\end{remark}

\begin{remark}
By the notion of $$
\d\xi(t) = \mathbf{v}^*(t,\xi(t))\d t + \sigma^*(t,\xi(t))\d W(t)
$$
we assume the probabilistic measure space $ \Omega $ and the motion space $ E $ are independent. Where the particles are does not affect how their stochastic parts behave. This is needed for the fact $$
\E\left(\int_t^{t+\delta(t)} \sigma^*(s,\xi(s))\d W(s)\bigg\vert \xi(t)\right) = \E\int_t^{t+\delta(t)} \sigma^*(s,\xi(s))\d W(s)
 = 0$$
to hold.
\end{remark}

\begin{remark}
Note that in our statement we have confined ourselves to the language of \eqref{StochasticProcess_Mainresults}, and this language requires the differentiability of the stochastic part of $ \xi $. We are proving effectively the fact that when the motion of $ \xi $ is separated into determined part and stochastic part, the variance of the stochastic part must be zero, and the stochastic part vanishes. Thus the determined part must be ODE \eqref{ODE}. The differentiability of the stochastic part is for the convenience of language but not essential for theories.
\end{remark}

\section{Proof of Main Results}
\textbf{In the first step of the proof} we prove particle motion \eqref{ODE} includes a statistical description \eqref{TheContinuityEquation}, the ODE is sufficient for the continuity PDE. In our proof, we consider an infinitely small time interval $ \Delta t $. It should be written as $ \Delta t $, but we choose to write differentiation $ dt $, indiscriminately representing a very small time interval considered. Consider location domain $ \Omega $ and time period $ (t,t+dt] $. Please note that $ \Omega $ is popular language usage in PDE theories to describe a spatial domain, and we shall later alter to the use of $ E $ to describe the location domain and reserve the symbol of $ \Omega $ for probabilistic measure space behind the considered stochastic motion.

$$ \int_\Omega \rho(x;t+dt)-\rho(x;t)\d x=\left(\int_\Omega \frac{\partial\rho}{\partial t} \d x\right) dt + o(dt)
$$ is the amount of "charge" increased in $ \Omega $ during time $ (t,t+dt] $, and this amount of charge is the amount of particles flowing in subtracting the amount of particles flowing out. $ o(dt) $ may be omitted from writing. Let $ \Omega_1 $ be the spatial location where particles flowing in $ \Omega $ are distributed at time $ t+dt $, and $ \Omega_2 $ be the spatial location where particles flowing out $ \Omega $ outside are distributed at time $ t+dt $. $ dt $ is an extremely short period such that particles counted either flow in or flow out. $ \Omega_1 $ and $ \Omega_2 $ do not intercept, and $ \Omega_1 + \Omega_2 $ together constructs an extremely narrow boundary of the considered location domain $ \Omega $. We have
\begin{align*}
\int_\Omega \rho(x;t+dt)-\rho(x;t)\d x=&\left(\int_\Omega \frac{\partial\rho}{\partial t} \d x\right)dt\\
=&\int_{\Omega_1}\rho(x;t+dt)\d x - \int_{\Omega_2}\rho(x;t+dt)\d x
\end{align*}
\includegraphics[width=0.95\textwidth]{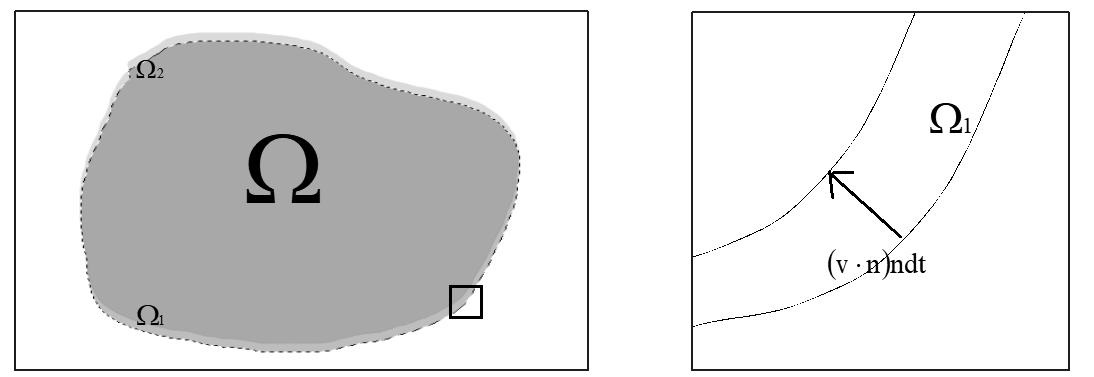}

Because $ dt $ is extremely short considering, $ \Omega_1 $ and $ \Omega_2 $ are very narrow boundary, and any locality of this boundary can be seen as rectangular infinitesimal element. We proceed the integration on a curved surface of little thickness in the $ p $-dimensional Euclidean space to the integration of curved surface in $ (p-1) $-dimensional space times an integration over the very short line segments. Mathematically it happens a substitution of variable. The problem that arises is that this substitution of variable is nowhere bijection when the surface curves. This mathematical problem can be fixed in two ways: One is to consider rectangular open to the left and closed to the right domains $ \Omega $ only, and the measure theory later holds if only rectangles are concerned. The other way is to strengthen constraints of the curvature of $ \partial\Omega $, the boundary of $ \Omega $, and yields an additional $ o(dt) $. We proceed
\begin{align}
&\int_{\Omega_1}\rho(x;t+dt)\d x - \int_{\Omega_2}\rho(x;t+dt)\d x \notag\\
=&\int_{\partial\Omega\cap\Omega_1} \left( \int_{y\in[s,s+(\mathbf{v}(t,s)\cdot\mathbf{n})\mathbf{n}dt]} \rho(y;t+dt)\d y \right)\d s \notag\\ 
&- \int_{\partial\Omega\cap\Omega_2} \left( \int_{y\in[s,s+(\mathbf{v}(t,s)\cdot\mathbf{n})\mathbf{n}dt]} \rho(y;t+dt)\d y \right)\d s \label{5y2}
\end{align}
Here $ s $ is a point on the curved surface $ \partial\Omega $, $ [s,s+(\mathbf{v}(t,s)\cdot\mathbf{n})\mathbf{n}dt] $ indicates variable $ y $ taking integration over one-dimensional line segment, and $ \mathbf{n} $ is the outward unit normal vector of the surface $ \partial\Omega $.

Assert \eqref{5y2} equals to \begin{align*}
=&\int_{\partial\Omega\cap\Omega_1} \left( \int_{y\in[s,s+(\mathbf{v}(t,s)\cdot\mathbf{n})\mathbf{n}dt]} \rho(s;t+dt)\d y +o(dt)\right)\d s \\
&- \int_{\partial\Omega\cap\Omega_2} \left( \int_{y\in[s,s+(\mathbf{v}(t,s)\cdot\mathbf{n})\mathbf{n}dt]} \rho(s;t+dt)\d y +o(dt)\right)\d s
\end{align*} 
This is because the line segment is $ O(dt) $ long, the integration over the line segment yields a $ O(dt^2) = o(dt) $. To show in strict mathematics, $$
\rho(y,t+dt) = \rho(x,t+dt) + \triangledown_x\rho(x+\xi (y-x),t+dt)\vert_{\xi\in[0,1]}\cdot(y-x)
$$and by integrating $ \triangledown_x\rho(x+\xi (y-x),t+dt)(y-x) $ yields
\begin{align*}
&\int_{y\in[x,x+(\mathbf{v}(t,x)\cdot\mathbf{n})\mathbf{n}dt]} \triangledown_x\rho(x+\xi (y-x),t+dt)\cdot(y-x)\d y \\
\leqslant& C \int_{y\in[x,x+(\mathbf{v}(t,x)\cdot\mathbf{n})\mathbf{n}dt]}(y-x)\d y \\
\leqslant & C(dt)^2
\end{align*}
where $ C $ is an absorbing constant. 

We proceed from \eqref{5y2} to:
\begin{align*}
=&\int_{\partial\Omega\cap\Omega_1} \left( \int_{y\in[x,x+(\mathbf{v}(t,x)\cdot\mathbf{n})\mathbf{n}dt]} \rho(x;t+dt)\d y +o(dt)\right)\d x \\
&- \int_{\partial\Omega\cap\Omega_2} \left( \int_{y\in[x,x+(\mathbf{v}(t,x)\cdot\mathbf{n})\mathbf{n}dt]} \rho(x;t+dt)\d y +o(dt)\right)\d x\\
=& \int_{\partial\Omega\cap\Omega_1 + \partial\Omega\cap\Omega_2} - \mathbf{n}\cdot \mathbf{v}(t,x)\rho(x;t+dt)\d xdt+o(dt)\\
=& \int_{\partial\Omega} - \mathbf{n}\cdot \mathbf{v}(t,x)\rho(x;t+dt)\d xdt+o(dt)\\
=&\int_\Omega -\mathrm{div}( \mathbf{v}(t,x)\rho(x;t+dt))\d xdt+o(dt)
\end{align*}
Notice at $ \partial\Omega\cap\Omega_1 $, $ \mathbf{v}\cdot\mathbf{n} $ is negative and there is additional negative sign yielded. Now we have proved that \begin{equation}\label{eq4.1.1}
\int_\Omega \frac{\partial\rho}{\partial t}\d x dt = \int_\Omega -\mathrm{div}( \mathbf{v}(t,x)\rho(x;t+dt))\d xdt+o(dt)
\end{equation}
by integrating \eqref{eq4.1.1} from time $ t_1 $ to time $ t_2 $ $$
\int_{[t_1,t_2]}\int_\Omega \frac{\partial\rho}{\partial t}\d x \d t = \int_{[t_1,t_2]}\int_\Omega -\mathrm{div}( \mathbf{v}(t,x)\rho(x;t))\d x\d t
$$
The last process is typical in measure theory, location-time domains like $ \Omega\times[t_1,t_2] $ constitutes semi-algebra in measure theory, and thus $$
\iint_A \frac{\partial\rho}{\partial t}\d x\d t = \iint_A -\mathrm{div}( \mathbf{v}(t,x)\rho(x,t))\d x\d t
$$for any location-time $ A $, which is equivalent to $$
\frac{\partial}{\partial t}\rho(x,t) + \mathrm{div}( \mathbf{v}(t,x)\rho(x,t)) = 0
$$This completes the proof of particle motion in velocity field containing its statistical nature; and we are proving the converse in all the text coming up next.

\textbf{In the second step of the proof} we prove the uniqueness of the time-dependent velocity field $ \mathbf{v}\in \mathbb{R}\times \mathbb{R}^p\to \mathbb{R}^p $ in the continuity equation. There cannot be two velocity fields satisfying the same one continuity equation. If $$
\frac{\partial}{\partial t}\rho(x;t) + \mathrm{div}\left(\rho(x;t) \mathbf{v}^1(t,x)\right)=0
$$ and $$
\frac{\partial}{\partial t}\rho(x;t) + \mathrm{div}\left(\rho(x;t) \mathbf{v}^2(t,x)\right)=0
$$ both hold for a statistical density $ \rho(x,t) $, then there must be $$
\mathbf{v}^1 \equiv \mathbf{v}^2
$$ 
To prove the result, we have\begin{equation}\label{Prf1.3.01}
\mathrm{div}\left((\mathbf{v}^1(t,x)-\mathbf{v}^2(t,x))\rho(x,t)\right)\equiv 0 
\end{equation}
Notice that time-dependent ODE \eqref{ODE} can define its time reversal, the continuity equation can have its time reversal$$
\frac{\partial}{\partial s}\rho(x,s) + \mathrm{div}( -\mathbf{v}(t-s,x)\rho(x,s)) = 0
$$ 
The fact that the continuity equation holds for any initial value $ \rho(x,0) $ equals to the fact that that holds for any of its time  $ \rho(x,t) $. Take $ \rho(x,t) = \mathrm{const} $ into \eqref{Prf1.3.01} and we have $$
\mathrm{div}\left((\mathbf{v}^1(t,x)-\mathbf{v}^2(t,x))\right)\equiv 0 
$$ Also from \eqref{Prf1.3.01} we have $$
\mathrm{div}\left(\mathbf{v}^1(t,x)-\mathbf{v}^2(t,x)\right)\rho(x,t) + (\mathbf{v}^1(t,x)-\mathbf{v}^2(t,x))\cdot \triangledown \rho(x,t)\equiv 0 
$$ which reduces to $$
(\mathbf{v}^1(t,x)-\mathbf{v}^2(t,x))\cdot \triangledown \rho(x,t)\equiv 0
$$ 
Take $ \rho(x,t) $ as specific form $ \rho(x,t) \equiv \rho_j(x_j), j=1,\cdots, p $ $$
\left(\mathbf{v}^1(t,x)-\mathbf{v}^2(t,x)\right)_j \frac{\d}{\d x_j}\rho_j(x_j)\equiv 0
$$ here $ \left(\mathbf{v}^1(t,x)-\mathbf{v}^2(t,x)\right)_j $ denotes the $ j $-th component of $ \left(\mathbf{v}^1(t,x)-\mathbf{v}^2(t,x)\right) $. By taking non-constant $ \rho_j(x_j) $, $$
\left(\mathbf{v}^1(t,x)-\mathbf{v}^2(t,x)\right)_j \equiv 0
$$ for any $ j $, and this completes the proof of $ \mathbf{v}^1 \equiv \mathbf{v}^2 $.

Though simple in mathematics, the second part of the proof is core to the whole theory. We are going to discuss the containing of particle nature by the continuity equation; however, if the velocity field by the continuity equation is not unique, there is no way we can discuss the containing of ODE by the continuity equation.

\textbf{In the third step of the proof} we are to study how a small iron-like air ball behaves in the air stream in a room. Note that this cannot be done by unelaborately substituting Dirac delta function into the continuity equation, because Dirac delta function violates the differentiability of mathematics. We borrow the probabilistic tools of expectation and variance. The classic Chebyshev inequality states that when a probabilistic distribution has variance of nearly zero, all of it would cumulate at its expectation. In a violation of correctness of probabilistic language, we are going to prove that the continuity equation \eqref{TheContinuityEquation} not only guarantees the uniqueness of the velocity field but also truly designates non-stochastic direction of motion for the cumulated mass: If \begin{equation}\label{Prf1.3.02}
\frac{\partial}{\partial t}\rho(u,t) + \mathrm{div}_u\left(\mathbf{v}(t,u)\rho(u,t)\right)=0
\end{equation}
holds for any macro density $ \rho(u,t) $, then the statement \begin{equation}\label{Prf1.3.05}
\text{"Particles would move in the direction of the velocity field }\mathbf{v} \text{ for a very instant time"}
\end{equation}
holds. Of course this language is not acceptable in probability, and we will give precise probabilistic languages in the fourth and fifth steps of the proof.

We use the symbol $ u $ as integration variable and $ x $ as the expectation of $ \rho(u;t) $ $$
\int u\rho(u;t)\d u = x
$$
We use the symbol $ \Sigma $ as the variance of $ \rho(u;t) $ $$
\int (u-x)(u-x)^T\rho(u;t)\d u = \Sigma
$$
We are to investigate the expectation and variance of $ \rho(u;t+dt) $ for a very instant time $ dt $ forward, where $$
\frac{\partial}{\partial t}\rho(u,t) + \mathrm{div}_u\left(\mathbf{v}(t,u)\rho(u,t)\right)=0
$$
the continuity equation is satisfied. Our result is that the expectation of $ \rho(\cdot,t+dt) $ is along the direction of the velocity $ \v(x,t) $ of $x $ with variance being arbitrarily small, when $ \Sigma $ the variance of $ \rho(\cdot;t) $, the previous-$dt$ state of $ \rho(\cdot;t+dt) $, is set arbitrarily small.

Consider the scenario that $ \Sigma $ is being set very small $ \Sigma\approx 0 $, and in this case $ \rho(\cdot;t) $ is almost a Dirac delta function which concentrates at $ x $. Typically we can consider $$
\Sigma = \sigma^2 \Sigma^{(0)} \text{ with } \sigma\approx 0
$$ $ \Sigma^{(0)} $ being a given fixed positive-definite matrix. We write $$
o(\Sigma) = O(\sigma^2) = o(\sigma)
$$
to omit the auxiliary symbol of $ \sigma $. 

In the text that follows, $ o(\Sigma) $, essentially being $ o(\sigma) $ if by assumption, shows how small the yield-out is by setting the variance of $ \rho(\cdot;t) $ $ \Sigma $, but not necessarily presumes the existence of limiting process. What we mean by $ o(\Sigma) $ is an item that eventually goes to zero by setting $ \Sigma $ the variance of $ \rho(\cdot;t) $ arbitrarily small. The readers may regard it happening a limiting process by considering sequentially many more and more shrinking $ \rho(\cdot;t) $ to Dirac delta function, and observe the shrinking phenomenon of $ \rho(\cdot;t) $ and $ \rho(\cdot;t+dt) $ . To present the text we do not care or arduously write out about limiting. 

We approach the target of \eqref{Prf1.3.05} in the following two lemmas: Lemma \ref{Prf1.3.03} and Lemma \ref{Prf1.3.06}.

\begin{lemma}\label{Prf1.3.03}
The expectation of $ \rho(\cdot;t+dt) $ is $$
\int u\rho(u;t+dt)\d u = x + \mathbf{v}(t,x)dt + o(\Sigma)dt + o(dt)
$$
\end{lemma}
\begin{prof}
The proof of Lemma \ref{Prf1.3.03} is as follows: Rewrite \eqref{Prf1.3.02} in the differentiation form:\begin{equation}\label{Prf1.3.04}
\rho(u;t+dt) = \rho(u;t) - \mathrm{div}\left(\mathbf{v}(t,u)\rho(u,t)dt\right) + o(dt)
\end{equation} Take expectation of both sides of \eqref{Prf1.3.04}, and look up at the right hand side. $$
\int u\rho(u;t)du = x
$$ is provided, and by using the Newton-Leibniz theorem or the Gauss-Green theorem \begin{align*}
&-\int u \mathrm{div}_u\left(\mathbf{v}(t,u)\rho(u,t)dt\right)\d u  \\
=& \int \triangledown_u u\cdot \left(\mathbf{v}(t,u)\rho(u,t)dt\right)\d u\\
=& \int \left(\mathbf{v}(t,u)dt\right)\rho(u,t)\d u
\end{align*}
Here $ -\int u_j \mathrm{div}_u\left(\mathbf{v}(t,u)\rho(u,t)dt\right)\d u = \int \mathbf{v}_j(t,u)dt\rho(u,t)\d u $, $\triangledown u $ yields identity matrix. Notice here $ u $ and $ x $ are $ p$-dimensional, and that the integrable variable is $ \d u $, $ dt $ is the differentiation already done before. Assert \begin{equation}\label{Prf1.3.04.1}
\int \left(\mathbf{v}(t,u)dt\right)\rho(u,t)du = \mathbf{v}(t,x)dt + o(\Sigma)dt
\end{equation}
Once this assertion is proved, the lemma is completed.

To prove \eqref{Prf1.3.04.1} we have that \begin{align}\label{eqn3}
& \int \left(\mathbf{v}(t,u)dt\right)\rho(u,t)\d u - \mathbf{v}(t,x)dt=\\ \notag
& \int_{B(x,\epsilon)} \left(\mathbf{v}(t,u)dt - \mathbf{v}(t,x)dt\right)\rho(u,t)\d u +\int_{B(x,\epsilon)^c}  \left(\mathbf{v}(t,u)dt- \mathbf{v}(t,x)dt\right)\rho(u,t)\d u 
\end{align}
for any $ \epsilon $,
\begin{align*}
&\int_{B(x,\epsilon)^c}  \left(\mathbf{v}(t,u)dt- \mathbf{v}(t,x)dt\right)\rho(u,t)\d u \\ 
\leqslant & 2\sup(\mathbf{|v|})dt \int_{B(x,\epsilon)^c} \rho(u,t)\d u\\
\leqslant & 2\sup(\mathbf{|v|})dt \int_{B(x,\epsilon)^c} \frac{(u-x)^T(u-x)}{\epsilon^2}\rho(u,t)\d u\\
\leqslant & 2\sup(\mathbf{|v|})dt  \frac{\mathrm{trace}(\Sigma)}{\epsilon^2} = \frac{o(\sigma)}{\epsilon^2} = \frac{o(\Sigma)}{\epsilon^2} 
\end{align*}
This requires that $ \mathbf{v} $ is uniformly bounded, $ \sup\mathbf{|v|} <\infty $; however, we are not cautious specifically of the boundedness condition of velocity field $ \mathbf{v} $ due to the following two reasons: One is that in physics the experiments concentrate at one locality, and it is meaningless to consider the field far beyond. The field should go to zero at very far. Second is that in numerical studies we take the whole measure of $ \rho $ in a compact set, condition stronger than the tightness of probability measure. By setting $ \Sigma $ sufficiently small
$$
\int_{B(x,\epsilon)^c}  \left(\mathbf{v}(t,u)dt- \mathbf{v}(t,x)dt\right)\rho(u,t)\d u $$
is arbitrarily close to zero. 

By the continuous condition of $ \mathbf{v} $ $$
\int_{B(x,\epsilon)} \left(\mathbf{v}(t,u)dt - \mathbf{v}(t,x)dt\right)\rho(u,t)\d u = o(\epsilon)
$$ Thus we have $$
\int \left(\mathbf{v}(t,u)dt\right)\rho(u,t)\d u - \mathbf{v}(t,x)dt \leqslant o(\epsilon) + 2\sup(\mathbf{|v|})\frac{\mathrm{tr}(\Sigma)}{\epsilon^2}
$$
By setting arbitrarily small $ \epsilon $ and then arbitrarily small $ \sigma $, the result of $  o(\epsilon) + 2\sup(\mathbf{|v|})\frac{\mathrm{tr}(\Sigma)}{\epsilon^2} $ deduces to $ o(\sigma) $. Thus the lemma is proved, and we still write $ o(\Sigma) $ in place of $ o(\sigma) $ to simplify the notation and emphasise it is a result of the variance set to $ \rho(\cdot,t) $.

\end{prof}

\begin{lemma}\label{Prf1.3.06}
The variance of $ \rho(\cdot;t+dt) $ is given by $$
\int (u-x-  \mathbf{v}(t,x)dt )(u-x-  \mathbf{v}(t,x)dt )^T\rho(u;t+dt)\d u = \Sigma + o(\Sigma)dt + o(dt)
$$
The coefficient before $ dt $ the linear order of Taylor expansion of the variance of $ \rho(\cdot;t+dt) $ is a $ o(\Sigma) $.
\end{lemma}
\begin{prof}
We proceed by analysing\begin{equation}\label{lem5.2varicance}
\int (u-x -  \mathbf{v}(t,x)dt - o(\Sigma)dt - o(dt))(u-x -  \mathbf{v}(t,x)dt - o(\Sigma)dt - o(dt))^T\rho(u;t+dt)\d u
\end{equation}still rewrite $ \rho(u;t+dt) $ in the differentiation form $$
\rho(u;t+dt) = \rho(u;t) - \mathrm{div}\left(\mathbf{v}(t,u)\rho(u,t)dt\right) + o(dt)
$$This is a little disastrous in writing, but the analysis is clear:
\begin{itemize}
\item $ \int (u-x)(u-x)^T\rho(u;t)\d u $ yields $ \Sigma $.
\item The $ dt $ item $ \int (u-x)(\mathbf{v}(t,x) + o(\Sigma))^T\rho(u;t)\d u = \int (u-x)\rho(u;t)\d u (\mathbf{v}(t,x) + o(\Sigma))^T $ yields $ 0 $.
\item We need only to tackle the remaining $ dt $ term being $ -\int (u-x)(u-x)^T\mathrm{div}\left(\mathbf{v}(t,u)\rho(u,t)dt\right)\d u $, because the rest are all $ o(dt) $.
\item Consider the $ (i,j) $-th component of $ (u-x)(u-x)^T $ and consider its integration $ -\int (u-x)_i(u-x)_j\mathrm{div}\left(\mathbf{v}(t,u)\rho(u,t)dt\right)\d u $, which equals to$$
\int \left((u-x)_i \mathbf{v}(t,u)_j+ (u-x)_j \mathbf{v}(t,u)_i\right)(\rho(u,t)dt)\d u
$$
Apply the same technique in the last lemma:\begin{align*}
&\int (u-x)_i \mathbf{v}(t,u)_j\rho(u,t)\d u = \\
&\int_{B(x,\epsilon)} (u-x)_i \mathbf{v}(t,u)_j\rho(u,t)\d u+\int_{B(x,\epsilon)^c} (u-x)_i \mathbf{v}(t,u)_j\rho(u,t)\d u
\end{align*}
$ \int_{B(x,\epsilon)} (u-x)_i \mathbf{v}(t,u)_j\rho(u,t)\d u $ is  $ o(\epsilon) $, and\begin{align*}
&\left| \int_{B(x,\epsilon)^c} (u-x)_i \mathbf{v}(t,u)_j\rho(u,t)\d u \right|\\
\leqslant & \int_{B(x,\epsilon)^c}\left|(u-x)_i \mathbf{v}(t,u)_j\right|\rho(u,t)\d u \\
\leqslant & \sup(|\v|) \int_{B(x,\epsilon)^c}\left|(u-x)_i \right|\rho(u,t)\d u \\
\leqslant & \sup(|\v|) \frac{\int_{B(x,\epsilon)^c}(u-x)_i^2 \rho(u,t)\d u}{\epsilon}
\sim \frac{\sigma^2}{\epsilon}
\end{align*}
Thus $$
\int (u-x)_i \mathbf{v}(t,u)_j\rho(u,t)\d u \leqslant o(\epsilon) + \sup(|\v|) \frac{\sigma^2}{\epsilon}
$$ 
and the result is the same for $$
\int \left( (u-x)_j \mathbf{v}(t,u)_i\right)(\rho(u,t)dt)\d u
$$
By setting arbitrarily small $ \epsilon $ and then arbitrarily small $ \sigma $, the result of $ o(\epsilon) + \sup(|\v|) \frac{\sigma^2}{\epsilon} $ is $ o(\sigma) = o(\Sigma) $
\end{itemize}
The lemma is proved. 
\end{prof}

From Lemma \ref{Prf1.3.03} and Lemma \ref{Prf1.3.06}, we have shown that the trajectory of velocity field $ \mathbf{v}(t,x) $ is taken by the small iron-like ball, compressed of the air in the whole room by our brainstorming, in the consideration of $ dt $ the linear order of time expansion. By setting $ \Sigma $ arbitrarily small, the whole mass appears at time $ t+dt $ in the direction $ \mathbf{v}(t,x)dt $ of $ x $ its last cumulation. This result seems obvious, since the continuity equation describes particles, it definitely acts in accordance with the particle motion when the scenario "compressing all air in the room into a small iron ball" is being considered.

\textbf{In the fourth step of the proof} we conclude our result in the formal probabilistic language. We are going to prove that for the diffusion process $ \xi(t) $ \begin{equation}
\label{DiffusionProcess2}
\d\xi(t) = \mathbf{v}^*(t,\xi(t))\d t + \sigma^*(t,\xi(t))\d B(t)
\end{equation}
$ B_t $ being Brownian motion. If $ \rho(\cdot;t) $ the probability density of $ \xi(t) $ satisfies the continuity equation \eqref{TheContinuityEquation}; Or equivalently for any independent and identically distributed diffusion processes $ \xi^{(1)}(t),\xi^{(2)}(t),\cdot,\xi^{(N)}(t), N = 1000 $ mol \begin{equation}
\label{DiffusionProcess3}
\d\xi^{(j)}(t) = \mathbf{v}^*(t,\xi^{(j)}(t))\d t + \sigma^*(t,\xi^{(j)}(t))\d B^{(j)}(t), j = 1,\cdots,N
\end{equation}
$ B^{(j)}_t $ being Brownian motion. If $ \rho $ the statistical density of $ \xi^{(j)}(t),j=1,\cdots,N $ satisfies the continuity equation \eqref{TheContinuityEquation}. Then we have $$
\mathbf{v}^*(\cdot,\cdot) \equiv \mathbf{v}(\cdot,\cdot) 
$$and $$
\sigma^*(\cdot,\cdot) \equiv 0 
$$ $ \xi $ is the probabilistically indistinguishable process with the trajectory of ODE of \eqref{ODE}. 

This result is included in the third part of the proof. Brownian motion has variance growth proportional to time elapse, which contradicts with the observation of $ o(\Sigma) $ coefficient before $ dt $.

\textbf{In the fifth and last step of the proof} we generalise our result from diffusion processes \eqref{DiffusionProcess2} to all stochastic processes. Still we reserve the symbol $ x $ to represent the expectation, and use $ u $ as the integrated variable. In this fifth part of proof, we presume $ \v(\cdot,t) \in C^\infty $ and $ \rho(\cdot,t)$ is analytic. We presume $ \v^*(t,\cdot) \in C^\infty $ and the solution of ODE is analytic. In fact, the last illustration requires many conditions to hold. Here we are giving a complete proof for arbitrary stochastic process $ W(t) $ instead of Brownian motion $ B(t) $, and arbitrary $ p $-dimension besides the special case $ p = 1 $.

We still consider the case that the variance of $ \rho(\cdot;t) $ is set arbitrarily small $$
\Sigma = \sigma^2\Sigma^{(0)}, \sigma\downarrow 0
$$
and write $ o(\Sigma) $ instead of $ o(\sigma) $ to reduce symbols and to emphasise this is the result of $ \Sigma $.

\begin{proposition}\label{Prf1.3.07}
$$
\int f(u)\rho(u;t)\d u = f(x) + o(\Sigma)
$$
where $ x $ being the expectation of $ \rho(\cdot;t) $ $$
x = \int u \rho(u;t)\d u
$$
and the variance of $ \rho(\cdot;t) $ $$ \int (u-x)(u-x)^T \rho(u;t)\d u = \Sigma 
$$ 
is being set arbitrarily small. $ f $ in this proposition can be scalar-valued, vector-valued, or matrix-valued uniformly bounded function.
\end{proposition}

\begin{prof}
We have showcased this technique before in the proof of Lemma \ref{Prf1.3.03} and Lemma \ref{Prf1.3.06}, now we use it again. We have that  \begin{align*}
&\int f(u)\rho(u;t)\d u - f(x)\\
=&\int \left(f(u) - f(x) \right)\rho(u;t)\d u \\
=& \int_{B(x,\epsilon)} \left(f(u) - f(x)\right)\rho(u,t)\d u  +\int_{B(x,\epsilon)^c} \left(f(u) - f(x)\right)\rho(u,t)\d u
\end{align*}
Apply mean value theorem \begin{align*}
\left|\int_{B(x,\epsilon)} \left(f(u) - f(x)\right)\rho(u,t)\d u\right| = &\left| \int_{B(x,\epsilon)} \left(\triangledown f(\eta)\cdot(u-x)\right)\rho(u,t)\d u\right|\\
\leqslant& \int_{B(x,\epsilon)} \left|\triangledown f(\eta)\cdot(u-x)\right|\rho(u,t)\d u\\
\leqslant& \int_{B(x,\epsilon)} |\triangledown f(\eta)|\cdot|(u-x)|\rho(u,t)\d u\\
\leqslant& \int_{B(x,\epsilon)} \max_{\eta\in B(x,\epsilon)}(|\triangledown f(\eta)|)\epsilon\rho(u,t)\d u\\
= & O(\epsilon)
\end{align*}
$ \int_{B(x,\epsilon)} \left(f(u) - f(x)\right)\rho(u,t)\d u $ is always a $ O(\epsilon) $, the constant $ C_1 $ before $ \epsilon $ is the maximum of the norm of the gradient of $ f $ near $ x $ $$
C_1 = \max_{\eta\in B(x,\epsilon)}(|\triangledown f(\eta)|)
$$ which is a constant depending only on $ f $ itself. 

The result $$ \int_{B(x,\epsilon)^c} \left(f(u) - f(x)\right)\rho(u,t)\d u $$ is $
\frac{o(\Sigma)}{\epsilon^2} $, which requires uniform boundedness condition, such as $f $ to be zero outside a compact set $ E \subset \R^p $:
\begin{align*}
&\left|\int_{B(x,\epsilon)^c} (f(u) - f(x))\rho(u,t)du\right|\\
\leqslant& \int_{B(x,\epsilon)^c} \left|(f(u) - f(x))\right|\rho(u,t)du \\
\leqslant& C_2  \int_{B(x,\epsilon)^c} \rho(u,t)du \\
\leqslant& C_2  \int_{B(x,\epsilon)^c} \frac{(u-x)^T(u-x)}{\epsilon^2}\rho(u,t)du\\
\leqslant& C_2 \frac{\mathrm{trace}\Sigma}{\epsilon^2}
\end{align*}
where $$
C_2 = 2 \max_{x\in \R^P} |f(x)|
$$

Now we have that $$
\int (f(u)-f(x))\rho(u,t)du = C_1\epsilon + C_2\frac{\mathrm{trace}\Sigma}{\epsilon^2}
$$ 
For any $ \epsilon*>0 $, take $$
\epsilon = \frac{\epsilon^*}{2C_1}
$$
and for this $ \epsilon $ take $ \Sigma $ small enough such that $$
C_2\frac{\mathrm{trace}\Sigma}{\epsilon^2} = \frac{\epsilon^*}{2}
$$
we complete showing the result that $$
\int f(u)\rho(u;t)\d u = f(x) + o(\Sigma)
$$
\end{prof}

We are studying the variance of $ \rho(\cdot,t+dt) $, by tackling the higher orders of Taylor expansion $$ 
\rho(x;t+dt) = \rho(x;t) + \frac{\partial}{\partial t}\rho(x;t)dt + \frac{1}{2}\frac{\partial^2}{\partial t^2}\rho(x;t)dt^2 + \frac{1}{6}\frac{\partial^3}{\partial t^3}\rho(x;t)dt^3 + \cdots
$$
to pursue our result. We have proved the coefficient before $ dt $ a $ o(\Sigma) $, and now the higher order items of $ dt^j,j\geqslant 2 $ of the variance of $ \rho(\cdot;t+dt) $ are coming into consideration. By simply calculating the derivatives, we have \begin{align*}
&\frac{\partial}{\partial t}\rho(u;t) = - \mathrm{div}\left(\rho(u;t) \mathbf{v}(u,t)\right)\\
&\frac{\partial^2}{\partial t^2}\rho(u;t) = -\mathrm{div}\left(\frac{\partial}{\partial t}\rho(u;t) \mathbf{v}(u,t) + \rho(u;t) \frac{\partial}{\partial t} \mathbf{v}(u,t)\right)\\
&\cdots\\
&\frac{\partial^n}{\partial t^n}\rho(u;t) = - \mathrm{div}\left(\sum\limits_{i = 0}^{n-1}\begin{pmatrix}
n-1\\
i
\end{pmatrix}\frac{\partial^i}{\partial t^i} \rho(u;t) \frac{\partial^{n-1-i}}{\partial t^{n-1-i}} \mathbf{v}(u,t)\right)\\
&\cdots
\end{align*}
$ \rho $ has arbitrary order of derivative with respect to $ t $ if $ \v $ has arbitrary order of derivative to $ t $. We need further that $ \rho $ is analytic to $ t $, which means $ \rho $ has power series with respect to $ t $ for a very instant time period forward. Lemma \ref{Prf1.3.06} has calculated the variance of $ \rho(\cdot;t+dt) $ of its linear item of $ dt $; when calculating higher orders of the variance of $ \rho(\cdot;t+dt) $, we need the following result:

\begin{proposition}\label{Prf1.3.08}
For any scalar-valued function $ f $, $$ 
\int f(u)\frac{\partial^j}{\partial t^j}\rho(u;t)\d u = \int \left(\left(\frac{\partial}{\partial t}+\v(u,t)\cdot\frac{\partial}{\partial u} \right)^j f(u)\right)\rho(u;t)\d u,j\geqslant 1
$$
Define operator $ D = (\partial_t + \v(u,t)\cdot\partial_u) = (\partial_t + \v(u,t)\cdot\triangledown)$, we write simply $$
\int f(u)\partial_t^j\rho(u;t)\d u = \int D^jf(u)\rho(u;t)\d u = \int D^{j-1}(\v(u,t)\cdot\triangledown f(u))\rho(u;t)\d u
$$
\end{proposition}

\begin{prof}
We prove by method of induction. By the Gauss-Green theorem, we have that 
\begin{align*}
&\int f(u)\frac{\partial}{\partial t}\rho(u;t)\d u\\
=&\int f(u)(-\mathrm{div}(\v(u,t)\rho(u,t))\d u \\
=&\int \v(u,t)\cdot\triangledown f(u)\rho(u,t)\d u\\
=&\int \left(\left(\frac{\partial}{\partial t}+\v(u,t)\cdot\frac{\partial}{\partial u} \right) f(u)\right)\rho(u;t)\d u
\end{align*}
Note that the $ \cdot $ between $ \v $ and $ \partial_u $ may be omitted from writing to save space. Induction Hypothesis: Assume $$ 
\int f(u)\frac{\partial^j}{\partial t^j}\rho(u;t)\d u = \int \left(\left(\frac{\partial}{\partial t}+\v(u,t)\frac{\partial}{\partial u} \right)^j f(u)\right)\rho(u;t)\d u
$$ is true for positive integer $ j $, then \begin{align*}
&\int f(u)\frac{\partial^{j+1}}{\partial t^{j+1}}\rho(u;t)\d u\\
=&\int f(u)\frac{\partial}{\partial t}\frac{\partial^{j}}{\partial t^{j}}\rho(u;t)\d u\\
=&\frac{\d}{\d t}\int f(u)\frac{\partial^{j}}{\partial t^{j}}\rho(u;t)\d u\\
=&\frac{\d}{\d t}\int \left(\left(\frac{\partial}{\partial t}+\v(u,t)\frac{\partial}{\partial u} \right)^j f(u)\right)\rho(u;t)\d u\\
=&\int \partial_t\left(\left(\frac{\partial}{\partial t}+\v(u,t)\frac{\partial}{\partial u} \right)^j f(u)\right)\rho(u;t)\d u + \int \left(\left(\frac{\partial}{\partial  t}+\v(u,t)\frac{\partial}{\partial u} \right)^j f(u)\right)\partial_t\rho(u;t)\d u\\
=&\int \partial_t\left(\left(\frac{\partial}{\partial t}+\v(u,t)\frac{\partial}{\partial u} \right)^j f(u)\right)\rho(u;t)\d u + \\
&\int \left(\left(\frac{\partial}{\partial t}+\v(u,t)\frac{\partial}{\partial u} \right)^j f(u)\right)\left(-\mathrm{div}\left( \v(u,t)\rho(u;t) \right)\right)\d u\\
=&\int (\partial_t+\v(u,t)\cdot\partial_u) \left(\left(\frac{\partial}{\partial t}+\v(u,t)\frac{\partial}{\partial u} \right)^j f(u)\right)\rho(u;t)\d u\\
=&\int \left(\left(\frac{\partial}{\partial t}+\v(u,t)\frac{\partial}{\partial u} \right)^{j+1} f(u)\right)\rho(u;t)\d u
\end{align*}
This concludes that the induction hypothesis is correct for any positive integer, thus the proof of proposition is complete.
\end{prof}

Here we introduce an operator $ D = (\partial_t + \v(t,x)\cdot\partial_x) $, or $ D = (\partial_t + \v(t,u)\cdot\partial_u) $. This is one operator $$
D = (\partial_t + \v(t,\circ)\cdot\partial_\circ)
$$
and $ \circ $ depends on which spatial variable is being used. When we exert this operator to a vector $ \begin{pmatrix}
a_1\\
a_2\\
\vdots\\
a_p
\end{pmatrix}$
, we mean $$
 (\partial_t + \v(t,x)\partial_x) \begin{pmatrix}
a_1\\
a_2\\
\vdots\\
a_p
\end{pmatrix} = \begin{pmatrix}
 (\partial_t + \v(t,x)\cdot\partial_x)a_1\\
 (\partial_t + \v(t,x)\cdot\partial_x)a_2\\
\vdots\\
 (\partial_t + \v(t,x)\cdot\partial_x)a_p
\end{pmatrix}
$$
$  (\partial_t + \v(t,x)\partial_x) $ is exerted on every single item of the vector, and there is a central dot between  $\v(t,x)$ and $\partial_x$. Before and later we might omit the central dot in  $ (\partial_t + \v(t,x)\cdot\partial_x) $ and simply write $ (\partial_t + \v(t,x)\partial_x) $ to save space. When we exert the symbol to matrix-valued function, it is the same that we exert the operator to every item of the matrix.

With the symbol of $ (\partial_t + \v(t,x)\cdot\partial_x) $ or write $ (\partial_t + \v(t,x)\cdot\triangledown) $, we can represent the higher orders of the expansion of the solution of the particle motion $$
\left\{\begin{aligned}
&\frac{\d}{\d t}x(t) = \v(t,x(t))\\
&x(t) = x
\end{aligned} \right.
$$
The solution is introduced by a symbol of $$
x(t+s) = x + g(x,t;s)
$$
where $$ g(x,t;s) =  \v(t,x)s + \sum\limits_{j=1}^\infty \frac{s^{j+1}}{(j+1)!} (\partial_t + \v(t,x)\partial_x)^j\v(t,x) 
$$
is the power series of the shift of the ordinary differential equation from time $ t $ to time $ t+s $.

\def\gx{g(x,t;dt)}
\def\gu{g(u,t;dt)}

Now we have defined the symbol $ g $: \begin{align*}
\gu =& \sum\limits_{j=1}^\infty \left(\left(\frac{\partial}{\partial t}+\v(u,t)\frac{\partial}{\partial u} \right)^j u\right) \frac{dt^j}{j!} \\
=& \v(u,t)dt + \sum\limits_{j=2}^\infty \left(\left(\frac{\partial}{\partial t}+\v(u,t)\frac{\partial}{\partial u} \right)^{j-1} \v(u,t)\right) \frac{dt^j}{j!}
\end{align*}
and
\begin{align*}
\gx =& \sum\limits_{j=1}^\infty \left(\left(\frac{\partial}{\partial t}+\v(t,x)\frac{\partial}{\partial x} \right)^j x\right) \frac{dt^j}{j!} \\
=& \v(t,x)dt + \sum\limits_{j=2}^\infty \left(\left(\frac{\partial}{\partial t}+\v(t,x)\frac{\partial}{\partial x} \right)^{j-1} \v(t,x)\right) \frac{dt^j}{j!}
\end{align*}
Note that $ \v(x,t) $ and $ \v(t,x) $ are exactly the same but different styles of language. The symbol of $ g $ is exhaustedly used in later text. $ g(\circ,t;dt) $ is the spatial shift from time $ t $ to time $ t+dt $ made by ODE $$\left\{\begin{aligned}
&\frac{\d}{\d t}x(t) = \v(t,x(t))\\
&x(t) = \circ
\end{aligned}\right. 
$$

Besides $ g $, another symbol we introduce is $ \delta(t) $: $ [t,t+\delta(t)] $ is the instant time period forward such that $ \rho $ and the solution of ODE are analytic, i.e. the power series converges. It is safe to write $
g(\circ,t;s)
$ for $
0\leqslant s\leqslant \delta(t)
$.

\begin{lemma}\label{Prf1.3.09E}
The expectation of $ \rho(\cdot;t+dt) $ is given by $$
\int u \rho(u;t+dt)\d u = x + \v(x,t)dt + \sum\limits_{j=1}^\infty \frac{dt^{j+1}}{(j+1)!} (\partial_t + \v(t,x)\partial_x)^j\v(x,t) + o(\Sigma)
$$
Since we have developed the symbol $$ \gx =  \v(t,x)dt + \sum\limits_{j=1}^\infty \frac{dt^{j+1}}{(j+1)!} (\partial_t + \v(t,x)\partial_x)^j\v(t,x) 
$$
the result can write as $$
\int u \rho(u;t+dt)\d u = x + \gx + o(\Sigma)
$$
\end{lemma}

\begin{prof}
By Proposition \ref{Prf1.3.08} 
\begin{align*}
&\int u\rho(u,t+dt)\d u\\
=&\int u\rho(u,t)\d u + \int u \sum\limits_{j=1}^\infty \frac{dt^j}{j!} \partial_t^j \rho(u,t)\d u\\
=& x + \int \sum\limits_{j=1}^\infty \left(\left(\frac{\partial}{\partial t}+\v(u,t)\frac{\partial}{\partial u} \right)^j u\right) \rho(u,t)\d u\frac{dt^j}{j!}\\
=& x + \int \gu \rho(u,t)\d u
\end{align*}
and by Proposition \ref{Prf1.3.07} $$
x + \int \gu \rho(u,t)\d u =  x + \gx + o(\Sigma)
$$
\end{prof}

Let us calculate the variance of $ \rho(\cdot;t+dt) $.

\begin{lemma}\label{Prf1.3.10var}
The variance of $ \rho(\cdot;t+dt) $ is given by \begin{align*}
&\int (u-x-\gx)(u-x-\gx)^T \rho(u;t+dt)\d u \\
=& \Sigma + o(\Sigma) + \int ((u-x)\gu^T + \gu (u-x)^T)\rho(u;t)\d u
\end{align*}
\end{lemma}

\begin{prof}
In fact the variance of $ \rho(\cdot;t+dt) $ should write as $$
\int (u-x-\gx-o(\Sigma))(u-x-\gx-o(\Sigma))^T \rho(u;t+dt)\d u 
$$
and this equals to 
$$
\int (u-x-\gx)(u-x-\gx)^T \rho(u;t+dt)\d u + o(\Sigma)\\
$$
Note that we omit $ o(\Sigma) $ possibly in the whole text. We have that 
\begin{align*}
&\int (u-x-\gx)(u-x-\gx)^T \rho(u;t+dt)\d u\\
=&\int (u-x)(u-x)^T\rho(u;t+dt)\d u  - \int ((u-x)\gx^T  + \gx(u-x)^T)\rho(u;t+dt)\d u \\
&+ \gx\gx^T
\end{align*}
By the last lemma on the expectation of $ \rho(\cdot;t+dt) $, $$
\int (u-x)\rho(u;t+dt)\d u  = \gx + o(\Sigma)
$$
and$$
\int (u-x)^T\rho(u;t+dt)\d u  = \gx^T + o(\Sigma)
$$
hold. We proceed
\begin{align*}
=&\int (u-x)(u-x)^T\rho(u;t+dt)\d u  - 2\gx\gx^T + \gx\gx^T + o(\Sigma)\\
=&\int (u-x)(u-x)^T\rho(u;t+dt)\d u  - \gx\gx^T 
\end{align*}
and $ o(\Sigma) $ is sometimes omitted.

By far we have shown \begin{align*}
&\int (u-x-\gx-o(\Sigma))(u-x-\gx-o(\Sigma))^T \rho(u;t+dt)\d u\\
=& \int (u-x)(u-x)^T\rho(u;t+dt)\d u - \gx\gx^T + o(\Sigma)
\end{align*}
and the rest of work is to calculate $ \int (u-x)(u-x)^T\rho(u;t+dt)\d u $:\begin{align*}
&\int (u-x)(u-x)^T \rho(u;t+dt)\d u \\
=&\int (u-x)(u-x)^T \sum\limits_{j=0}^\infty \frac{dt^j}{j!} \frac{\partial ^j}{\partial t^j} \rho(u,t)\d u\\
=&\Sigma +\int (u-x)(u-x)^T \sum\limits_{j=1}^\infty \frac{dt^j}{j!} \frac{\partial ^j}{\partial t^j} \rho(u,t)\d u\\
=&\Sigma +\int \sum\limits_{j=1}^\infty \left(\frac{\partial}{\partial t}+\v(t,u)\frac{\partial}{\partial u} \right)^j\left((u-x)(u-x)^T \right)\frac{dt^j}{j!}\rho(u,t)\d u
\end{align*}
By doing binomial expansion \begin{align*}
&\int \sum\limits_{j=1}^\infty \left(\left(\frac{\partial}{\partial t}+\v(t,u)\frac{\partial}{\partial u} \right)^j((u-x)(u-x)^T  \right)\frac{dt^j}{j!} \rho(u,t)\d u\\
=&\int \sum\limits_{j=1}^\infty \left(\left(\frac{\partial}{\partial t}+\v(t,u)\frac{\partial}{\partial u} \right)^j(u-x)\right)\frac{dt^j}{j!}(u-x)^T  \rho(u,t)\d u \\
+&\int (u-x)\left(\sum\limits_{j=1}^\infty \left(\frac{\partial}{\partial t}+\v(t,u)\frac{\partial}{\partial u} \right)^j(u-x)^T  \right)\frac{dt^j}{j!}\rho(u,t)\d u \\
+&\int \sum\limits_{j=2}^\infty \sum\limits_{k=1}^{j-1}  \begin{pmatrix}
j\\
k
\end{pmatrix}\left((\partial_t+\v(t,u)\partial_u)^k(u-x)\right)\left((\partial_t+\v(t,u)\partial_u)^{j-k}(u-x)^T\right)\frac{dt^j}{j!} \rho(u,t)\d u\\
=& \int\left( (u-x)\gu^T + \gu(u-x)^T\right) \rho(u;t)\d u\\
+&\int \sum\limits_{j=2}^\infty \sum\limits_{k=1}^{j-1}  \left(\frac{dt^k}{k!} (\partial_t+\v(t,u)\partial_u)^k(u-x)\right)\left(\frac{dt^{j-k}}{(j-k)!} (\partial_t+\v(t,u)\partial_u)^{j-k}(u-x)^T\right)\rho(u,t)\d u
\end{align*}

\begin{remark}
The last manipulation is \begin{align*}
& \sum\limits_{j=1}^\infty \left(\frac{\partial}{\partial t}+\v(t,u)\frac{\partial}{\partial u} \right)^j (AB)\frac{dt^j}{j!} \\
=&B\left(\frac{\partial}{\partial t}+\v(t,u)\frac{\partial}{\partial u} \right)Adt + A\left(\frac{\partial}{\partial t}+\v(t,u)\frac{\partial}{\partial u} \right)Bdt\\
&+  \sum\limits_{j=2}^\infty \left(\frac{\partial}{\partial t}+\v(t,u)\frac{\partial}{\partial u} \right)^j (AB)\frac{dt^j}{j!}\\
=&B\left(\frac{\partial}{\partial t}+\v(t,u)\frac{\partial}{\partial u} \right)Adt + A\left(\frac{\partial}{\partial t}+\v(t,u)\frac{\partial}{\partial u} \right)Bdt\\
&+  \sum\limits_{j=2}^\infty \sum\limits_{k=0}^{j}  \left(\frac{dt^k}{k!} (\partial_t+\v(t,u)\partial_u)^kA\right)\left(\frac{dt^{j-k}}{(j-k)!} (\partial_t+\v(t,u)\partial_u)^{j-k}B)\right)\\
=&B \sum\limits_{j=1}^\infty \left(\frac{\partial}{\partial t}+\v(t,u)\frac{\partial}{\partial u} \right)^j (A)\frac{dt^j}{j!} + A \sum\limits_{j=1}^\infty \left(\frac{\partial}{\partial t}+\v(t,u)\frac{\partial}{\partial u} \right)^j (B)\frac{dt^j}{j!}\\
&+  \sum\limits_{j=2}^\infty \sum\limits_{k=1}^{j-1}  \left(\frac{dt^k}{k!} (\partial_t+\v(t,u)\partial_u)^kA\right)\left(\frac{dt^{j-k}}{(j-k)!} (\partial_t+\v(t,u)\partial_u)^{j-k}B)\right)
\end{align*}
and we change $$
\sum\limits_{j=2}^\infty \sum\limits_{k=1}^{j-1} = \sum\limits_{k=1}^\infty \sum\limits_{j=k+1}^{\infty} 
$$
It would be clear to see for $ A $ and $ B $ being scalars, but in fact $ A $ and $ B $ here are vectors. \qedsymbol
\end{remark}

We proceed
\begin{align*}
=& \int\left( (u-x)\gu^T + \gu(u-x)^T\right) \rho(u;t)\d u\\
+&\int \sum\limits_{k=1}^\infty \sum\limits_{j=k+1}^\infty \left(\frac{dt^k}{k!} (\partial_t+\v(t,u)\cdot\partial_u)^k(u-x)\right)\left(\frac{dt^{j-k}}{(j-k)!} (\partial_t+\v(t,u)\cdot\partial_u)^{j-k}(u-x)^T\right)\rho(u,t)\d u\\
=& \int\left( (u-x)\gu^T + \gu(u-x)^T\right) \rho(u;t)\d u + \int \gu\gu^T\rho(u;t)\d u
\end{align*}

Now we have reached our final result \begin{align*}
&\int (u-x-\gx-o(\Sigma))(u-x-\gx-o(\Sigma))^T \rho(u;t+dt)\d u\\
=& \int (u-x)(u-x)^T\rho(u;t+dt)\d u - \gx\gx^T + o(\Sigma)\\
=& \Sigma + \int\left( (u-x)\gu^T + \gu(u-x)^T\right) \rho(u;t)\d u\\
 &+ \int \gu\gu^T\rho(u;t)\d u - \gx\gx^T + o(\Sigma)
\end{align*}
Since the item $$
\int \gu\gu^T\rho(u;t)\d u - \gx\gx^T
$$
produces additional $ o(\Sigma) $, the lemma is proved.
\end{prof}

Now we have done all the preparation, and next is the kick:

\begin{proposition}
The statistical density $ \rho $ of i.i.d. stochastic processes $ \xi(t) $
\begin{equation}
\label{StochasticProcess}
\d\xi(t) = \mathbf{v}^*(t,\xi(t))\d t + \sigma^*(t,\xi(t))\d W(t)
\end{equation}
where $ W(t) $ is any stochastic process,
satisfies the continuity equation $$
\frac{\partial}{\partial t}\rho(u,t) + \mathrm{div}_u\left(\mathbf{v}(u,t)\rho(u,t)\right)=0
$$ for any of its possible initial values $ \rho(\cdot;0) $, where $ \rho(\cdot;0) $ is the probability density of $ \xi(0) $. Then we have $$
\sigma^*(\cdot,\cdot) \equiv 0 
$$
\end{proposition}

\begin{prof}
Review that $ \delta(t) $ is the time period that $ \rho(\cdot;t+s) $ is analytic, which is to have infinitely many Taylor's expansion items called power series, for $ 0\leqslant s\leqslant \delta(t) $. We have that \begin{align*}
\mathrm{var} (\xi(t + \delta(t))) = \E(\mathrm{var} (\xi(t+\delta(t))\vert \xi(t))) + \mathrm{var} (\E(\xi(t + \delta(t))\vert \xi(t)))
\end{align*}
We have
$$
\E(\xi(t + \delta(t))\vert \xi(t)) = \xi(t) + g^*(\xi(t),t;\delta(t))
$$\begin{align*}
g^*(x,t;s) = & \v^*(t,x)s + \sum\limits_{j=1}^\infty \frac{s^{j+1}}{(j+1)!} (\partial_t + \v^*(t,x)\partial_x)^j\v^*(t,x) \\
=& \sum\limits_{j=0}^\infty \frac{s^{j+1}}{(j+1)!} (\partial_t + \v^*(t,x)\partial_x)^j\v^*(t,x)
\end{align*}
and\begin{align*}
g^*(u,t;s) = & \v^*(t,u)s + \sum\limits_{j=1}^\infty \frac{s^{j+1}}{(j+1)!} (\partial_t + \v^*(t,u)\partial_u)^j\v^*(t,u) \\
=& \sum\limits_{j=0}^\infty \frac{s^{j+1}}{(j+1)!} (\partial_t + \v^*(t,u)\partial_u)^j\v^*(t,u)
\end{align*}

Let us calculate $ \mathrm{var} (\E(\xi(t + \delta( t))\vert \xi(t))) $. We have that \begin{align*}
&\mathrm{var} (\E(\xi(t + \delta( t))\vert \xi(t))) \\
=& \mathrm{var} ( \xi(t) + g^*(\xi(t),t;\delta( t)))  \\
=&\int (u + g^*(u,t;\delta(t)) - x - g^*(x,t;\delta(t))(u + g^*(u,t;\delta(t)) - x - g^*(x,t;\delta(t))^T \rho(u,t)\d u
\end{align*}
Because $$
\int g^*(u,t;\delta(t))\rho(u;t)\d u =  g^*(x,t;\delta(t))) + o(\Sigma)
$$
and $$
\int g^*(u,t;\delta(t))g^*(u,t;\delta(t))^T\rho(u;t)\d u =  g^*(x,t;\delta(t)))g^*(x,t;\delta(t)))^T + o(\Sigma)
$$
we have that
$$
\int g^*(u,t;\delta(t)) - g^*(x,t;\delta(t)))( g^*(u,t;\delta(t)) - g^*(x,t;\delta(t)))^T \rho(u,t)\d u = o(\Sigma)
$$
Notice $$
\int (u-x)g^*(x,t;\delta(t))^T\rho(u,t)\d u = 0
$$
Thus  \begin{align*}
&\mathrm{var} (\E(\xi(t + \delta( t))\vert \xi(t))) \\
=& \int (u-x)(u-x)^T \rho(u;t)\d u +\\
& \int (u-x)g^*(u,t;\delta(t))^T\rho(u;t)\d u + \int g^*(u,t;\delta(t))(u-x)^T\rho(u;t)\d u + o(\Sigma)\\
=& \Sigma + \int (u-x)g^*(u,t;\delta(t))^T\rho(u;t)\d u + \int g^*(u,t;\delta(t))(u-x)^T\rho(u;t)\d u + o(\Sigma)
\end{align*}

Now let us turn to $ \E (\mathrm{var}(\xi(t + \delta( t))\vert \xi(t))) $: Since the stochastic motion is independent with $ \xi(t)$
\begin{align*}
&\mathrm{var} (\xi(t+\delta(t))\vert \xi(t)) \\
&= \E\left(\left( \xi(t+\delta(t)) - \E(\xi(t+\delta(t))\vert \xi(t))\right) \left( \xi(t+\delta(t)) - \E(\xi(t+\delta(t))\vert \xi(t))\right)^T\big\vert \xi(t)\right) \\
& = \E\left(\left(\tilde{W} (t,t+\delta(t)) \right)\left(\tilde{W} (t,t+\delta(t))  \right)^T\vert\xi(t)\right)
\end{align*}
where $$
\tilde{W} (t,t+\delta(t)) = \int_t^{t+\delta(t)}\sigma^*(s,\xi(s))\d W(s)
$$
is the stochastic part of motion during $ [t,t+\delta(t)] $; and $$
\E(\mathrm{var} (\xi(t+\delta(t))\vert \xi(t))) = \var(\tilde{W}(t,t+\delta(t)))
$$
is its variance.

Suppose there exists $ x^* $ such that $$
\sigma^*(t,x^*) > 0
$$
by letting $ \xi(t) $ concentrate to that $ x^* $ we have that \begin{align*}
\var(\tilde{W}(t,t+\delta(t)))
\approx \sigma^*(t,x^*)\var(W(t+\delta(t)) - W(t))\sigma^*(t,x^*)^T
\end{align*}
However, by comparing the variance of $ \rho(\cdot;t+dt) $, given by Lemma \ref{Prf1.3.10var}, and the variance of $ \xi(t+dt) $, there are only items $ \Sigma $, $ o(\Sigma) $, $ \int (u-x)\gu^T\rho(u;t)\d u $, $ \int \gu (u-x)^T\rho(u;t)\d u $, $ \int (u-x)g^*(u,t;dt)^T\rho(u;t)\d u $, $ \int g^*(u,t;dt) (u-x)^T\rho(u;t)\d u $. All these items denies the existence of the item $$
\sigma^*(t,x^*)\var(W(t+\delta(t)) - W(t))\sigma^*(t,x^*)^T 
$$
which is not absorbed by $ \Sigma $. Thus we have denied the existence of $$
\sigma^*(t,x^*) > 0
$$
for its provision of alien variance.
\end{prof}

By the comparison of the variance of $ \rho(\cdot;t+dt) $ and the variance of $ \xi(t+\delta(t)) $, $$
\int (u-x-\gx)(u-x-\gx)^T \rho(u;t+dt)\d u
$$
and $$
\int (u + g^*(u,t;\delta(t)) - x - g^*(x,t;\delta(t))(u + g^*(u,t;\delta(t)) - x - g^*(x,t;\delta(t))^T \rho(u,t)\d u 
$$yield$$
 \Sigma + o(\Sigma) + \int \left((u-x)\gu^T + \gu (u-x)^T\right)\rho(u;t)\d u = o(\Sigma)
$$ and $$
\Sigma + o(\Sigma) + \int \left((u-x)g^*(u,t;\delta(t))^T + g^*(u,t;\delta(t)) (u-x)^T\right)\rho(u;t)\d u = o(\Sigma)
$$
To set $ dt = \delta(t)$, their difference $ o(\Sigma) $ eradicates the possibility of the existence of stochastic motion. Now the symbol we reserved for probability space come to use. Let $ \Omega $ be the probability space. There exists $ \Omega_t $ such that $$
\P(\Omega_t) = 1
$$
and particles do non-stochastic motion in the time period $[t,t+\delta(t)] $, a very instant time period ahead of $ t $ on $ \omega\in\Omega_t $. Finally there are countably many rational points in the time period $ [0, T]\subset \mathbb{R}^+ $, then there exists $$ \Omega^* = \bigcap_{t\in \mathbb{Q}\cap [0,T]}\Omega_t \subset \Omega, \P(\Omega^*) = 1 
$$ 
$ \mathbb{Q} $ stands for all rational numbers, and $ T $ can take $ \infty $, such that particles do determined motion in the whole time process for $ \omega\in\Omega^* $. We conclude that $ \xi $ is non-stochastic process almost surely. 

By the uniqueness of the velocity field by the second step of the proof,$$
\frac{\partial}{\partial t}\rho(u;t) + \triangledown\cdot\left(\rho(u;t) \mathbf{v}(u,t)\right)=0
$$
and $$
\frac{\partial}{\partial t}\rho(u;t) + \triangledown\cdot\left(\rho(u;t) \mathbf{v}^*(t,u)\right)=0
$$
both hold, and we have $$
\v\equiv \v^*
$$
and $ \xi $ is probabilistically indistinguishable with the ODE $$
\frac{\d x(t)}{\d t} = \v(t,x(t))
$$
By far we have proved the whole theorem.

\section{Discussion in Addition}
In the first part of the proof of main results, we choose a proof of Reynolds transport theorem different from the usual approach of conservation of mass in textbooks.

Proposition \ref{Prf1.3.07} is key to the proof of main results.  This result is easy to understand. By setting the variance of $ \rho(\cdot;t) $ sufficiently small, the whole mass is condensed to a neighbourhood of its center.

Proposition \ref{Prf1.3.08} is key to the proof of main results, but the meaning can be difficult to read. We require the knowledge of techniques of asymptotic analysis to give an understanding of this proposition. The motion of particles in time period $ [t,t+dt] $ is given by $$
(\circ)(t+dt) - (\circ)(t) = \sum\limits_{j=1}^\infty (\partial_t + \v(x,t)\cdot\partial_x)^j(\circ)\frac{dt^j}{j!}
$$
where $ (\circ) $ representing particle motion $ x(t) $, and the shift of the statistical density of particles in time period $ [t,t+dt]$ is given by $$
\rho(\cdot;t+dt)-\rho(\cdot;t) = \sum\limits_{j=1}^\infty (\partial_t)^j\rho(\cdot;t)\frac{dt^j}{j!}
$$
To use classic symbols in asymptotic analysis, we replace $ dt $ by $ \alpha $: The motion of particles is given by $$
\sum\limits_{j=1}^\infty (\partial_t + \v(x,t)\cdot\partial_x)^j(\circ)\frac{\alpha^j}{j!}
$$
and the shift of their statistical density is given by $$
\rho(\cdot;t+dt) = \sum\limits_{j=1}^\infty (\partial_t)^j\rho(\cdot;t)\frac{\alpha^j}{j!}
$$
According to general principles of asymptotic analysis, $ \alpha $ to the same order should be classified together: the $ k $-th order of particle motion $$
(\partial_t + \v(x,t)\cdot\partial_x)^k(\circ)\frac{\alpha^k}{k!}
$$
should introduce to the effect of the $ k $-th order of shift to $ \rho $: $$
(\partial_t)^k\rho(\cdot;t)\frac{\alpha^k}{k!}
$$
This idea to be stated out by strict mathematics is our proposition.


\end{document}